\documentclass[twoside,12pt]{cmft}

\pagespan{1}{?}
\cmftinfo{XXYYYZZ}

\newtheorem{theorem}{Theorem}
\newtheorem{lemma}[theorem]{Lemma}
\newtheorem{corollary}[theorem]{Corollary}

\theoremstyle{remark}

\newcommand{\ip}[2]{{\left<#1,#2\right>}}
\newcommand{\tK}{\widetilde{K}}
\newcommand{\bbB}{{\mathbb B}}
\newcommand{\bbC}{{\mathbb C}}
\newcommand{\bbD}{{\mathbb D}}
\newcommand{\cH}{{\mathcal H}}
\newcommand{\cK}{{\mathcal K}}
\newcommand{\cL}{{\mathcal L}}
\newcommand{\lhk}{{\cL(\cH,\cK)}}
\newcommand{\lkh}{{\cL(\cK,\cH)}}
\newcommand{\lh}{{\cL(\cH)}}
\newcommand{\lfg}{{\cL(\cF,\cG)}}
\newcommand{\vertspace}{\vphantom{X^{X^X}}}
\newcommand{\ipzz}{{\left<z,z\right>}}
\newcommand{\cF}{{\mathcal F}}
\newcommand{\cG}{{\mathcal G}}

\begin{document}

\title[Schwarz-Pick inequalities]{Schwarz-Pick inequalities for the 
Schur-Agler class \\ on the polydisk and unit ball}
\date{}

\author[Anderson]{J.\ Milne Anderson}
\email{helen@math.ucl.ac.uk}

\address{Department of Mathematics,
  University College London,
  Gower Street,
  London  WC1E 6BT,
  UK}

\author[Dritschel]{Michael A.\ Dritschel}
\email{m.a.dritschel@ncl.ac.uk}
\address{School of Mathematics and Statistics,
  Merz Court \\ University of Newcastle upon Tyne,
  Newcastle upon Tyne  NE1 7RU,
  UK}

\author[Rovnyak]{James Rovnyak}
\email{rovnyak@virginia.edu}
\address{Department of Mathematics,
  University of Virginia,
  P.\ O.~Box 400137,
  Charlottesville, VA 22904,
  USA}


\thanks{J.\ Milne Anderson acknowledges the support of the 
Leverhulme Trust.}
\dedicatory{To Walter Hayman on his $80$-th birthday}

\begin{abstract}
  The notion of a unitary realization is used to estimate derivatives
  of arbitrary order of functions in the Schur-Agler class on the
  polydisk and unit ball.
\end{abstract}

\maketitle

\section{Introduction}\label{S:introduction}

The classical Schwarz-Pick estimate is the inequality
\begin{equation*}\label{aug22a}
  |\varphi'(z)| \le \frac{1 - |\varphi(z)|^2}{1-|z|^2}, \qquad |z|<1,
\end{equation*}
for an analytic function $\varphi(z)$ which is bounded by one on the
unit disk of the complex plane.  Ruscheweyh \cite{ruscheweyh} has
obtained best-possible estimates of higher order derivatives of
bounded analytic functions on the disk.  Similar estimates were
derived by other methods and for different classes of analytic
functions in one and several variables by Anderson and Rovnyak
\cite{AR}, Avkhadiev and Wirths \cite{AW2003,AW2004,AW2006},
B{\'e}n{\'e}teau, Dahlner, and Khavinson \cite{BDK}, and MacCluer,
Stroethoff, and Zhao \cite{MSZ1,MSZ2}.  In this paper we derive
estimates of derivatives of functions on the polydisk and unit ball
using realizations of functions as transfer functions of systems or,
equivalently, as characteristic functions of operator colligations.

The notion of a realization plays a key role in many areas of operator
theory, linear systems theory, and interpolation problems, where the
class of analytic functions which are bounded by one on the unit disk
$\bbD = \{ z \colon |z|<1 \}$ is often called the {\bf Schur class}.
By a {\bf realization} of an analytic function $\varphi(z)$ in the
Schur class on $\bbD$ is meant a representation in the form
\begin{equation}
  \label{aug22b}
  \varphi(z) = D + z C (I_\cH - zA)^{-1} B, \qquad z \in \bbD,
\end{equation}
where $A,B,C,D$ are the entries of a block operator matrix
\begin{equation}
  \label{aug22c}
    U=\begin{pmatrix} A & B \\[2pt]  C & D \end{pmatrix} \colon
  \begin{matrix} \cH \\[-3pt] \oplus \\[-1pt] \bbC \end{matrix} \to 
  \begin{matrix} \cH \\[-3pt] \oplus \\[-1pt] \bbC \end{matrix}
\end{equation}
for some Hilbert space $\cH$.  Such a representation always exists and
can be chosen such that $U$ is a unitary operator on $\cH \oplus \bbC$
(for example, see Brodski{\u\i} \cite[Theorem 5.1]{MSB}).  In this
case \eqref{aug22b} is called a {\bf unitary realization} of
$\varphi(z)$.  It turns out that the estimates of Ruscheweyh
\cite{ruscheweyh} have straightforward proofs by means of
\eqref{aug22b}, and moreover this method extends to several variables.

The notion of a unitary realization can be introduced for a broad
class of domains $\Omega$ in~$\bbC^d$.  The {\bf Schur class} on
$\Omega$ is the set of analytic functions which are defined and
bounded by one on~$\Omega$.  A {\bf unitary realization} of a function
$\varphi(z)$ in the Schur class on $\Omega$ is a representation
\begin{equation}\label{may29a}
\varphi(z) = D + C Z(z) (I_\cK - A Z(z))^{-1} B , \qquad z\in\Omega,
\end{equation}
where
\begin{equation}\label{may30a}
      U=\begin{pmatrix} A & B \\[2pt]  C & D \end{pmatrix} \colon
  \begin{matrix} \cH \\[-3pt] \oplus \\[-1pt] \bbC \end{matrix} \to 
  \begin{matrix} \cK \\[-3pt] \oplus \\[-1pt] \bbC \end{matrix}
\end{equation}
is a unitary operator for some Hilbert spaces $\cH$ and $\cK$, and
$Z(z)$ is a function of $z = (z_1,\dots,z_d)$ in $\Omega$ with
values in $\lkh$ such that
\begin{enumerate}
\item[(i)] $\| Z(z) \| < 1$ for all $z\in\Omega$.
\end{enumerate}
For the domains that we consider, we also assume:
\begin{enumerate}
\item[(ii)] $Z(z_1,\dots,z_d) = z_1 E_1 + \cdots + z_d E_d$ where $E_j
  \in\lkh$ and $\| E_j \| \le 1$, $j=1,\dots,d$.
\end{enumerate}
Condition (i) assures that the inverse in \eqref{may29a} exists for
each $z\in\Omega$.  Condition (ii) is assumed for convenience: it is
satisfied in our applications, and it simplifies differentiation
formulas because higher order derivatives of the linear function
$Z(z)$ vanish.  The estimates that we derive below generalize to
nonlinear functions $Z(z)$, but the nonlinear case is more complicated
because the formulas for higher order derivatives involve additional
terms.  As this generality is not needed here, we only treat the
linear case.  

Unitary realizations of the form \eqref{may29a} include
\eqref{aug22b} as a special case when $d=1$ and $\Omega = \bbD$.
We describe the situation for the polydisk $\bbD^d = \bbD \times
\cdots \times \bbD$ and unit ball $\bbB_d = \{ (z_1,\dots,z_d) \colon
|z_1|^2 + \cdots + |z_d|^2 < 1 \}$ in greater detail.

\subsection{The Schur-Agler class on the polydisk}\label{polydiskcase}
Let $\Omega = \bbD^d$.  A complex-valued analytic function
$\varphi(z)$ on $\bbD^d$ is said to belong to the {\bf Schur-Agler
  class on the polydisk} if for every set $T_1,\dots,T_d$ of commuting
contractions on a Hilbert space and every positive number $r<1$, $\|
\varphi(rT_1,\dots,rT_d) \| \le 1$.  For $d=1$ and $d=2$, the
Schur-Agler class coincides with the Schur class on $\bbD^d$.  For
$d>2$, the Schur-Agler class is a proper subset of the Schur class on
$\bbD^d$ (for example, see \cite{lottosteger} and \cite{varopoulos}).
For any $d \ge 1$, the class of functions which have a unitary
realization \eqref{may29a} coincides with the Schur-Agler class on
$\bbD^d$ when suitable choices are made.  Namely, in \eqref{may29a}
and \eqref{may30a} we take $\cK=\cH$ and assume that
\begin{equation}\label{may30c}
\cH = \cH_1 \oplus \cdots \oplus \cH_d 
\end{equation}
for some closed subspaces $\cH_1,\dots,\cH_d$ of $\cH$.  We also
assume that
\begin{equation}\label{jun2a}
Z(z_1,\dots,z_d) = z_1 E_1 + \cdots + z_d E_d ,
\end{equation}
where $E_1,\dots,E_d$ are the orthogonal projections of $\cH$ onto
$\cH_1,\dots,\cH_d$.  Alternatively, we can write the elements $h$ of
$\cH$ in column form:
\begin{equation*}
h = \begin{pmatrix}
 h_1 \\ \vdots \\ h_d
\end{pmatrix} .
\end{equation*}
Then $A,B,C$ can correspondingly be written in block operator form,
and in this notation,
\begin{equation*}
Z(z_1,\dots,z_d)
=
\begin{pmatrix}
 z_1 I_{\cH_1} & 0 & \cdots & 0 \\
   0 & \!\!\!z_2 I_{\cH_2} & \cdots & 0 \\
    & &\cdots & \\
   0 & 0 & \cdots & \!\!\!z_d I_{\cH_d}
\end{pmatrix} .
\end{equation*}
With these conventions, every function of the form \eqref{may29a}
belongs to the Schur-Agler class on $\bbD^d$, and conversely every
function in the Schur-Agler class on $\bbD^d$ has this form.  See
Agler \cite{agler1990} and Ball and Trent \cite{balltrent}.

\subsubsection*{Remark} The Schur-Agler class occurs naturally in the
function theory of the polydisk.  It is known that for every
function $\varphi(z)$ in the Schur class on $\bbD^d$,
\begin{equation*}  
(1 - |z_1|^2) \left| \frac{\partial \varphi}{\partial z_1} \right|
+ \cdots
+   (1 - |z_d|^2) \left| \frac{\partial \varphi}{\partial z_d} \right|
\le
1 - |\varphi(z)|^2 
\end{equation*}
at every point $z = (z_1,\dots,z_d)$ of $\bbD^d$ (see Rudin \cite[p.\ 
179]{rudinpolydisk} and Knese \cite{knese}).  The cases of equality are
determined by Knese \cite{knese}.  All extremal functions belong to
the Schur-Agler class.  A Schur function $\varphi(z)$ on $\bbD^d$,
which has no identically vanishing first partials, satisfies
\begin{equation*}
  (1 - |z_1|^2) \left| \frac{\partial \varphi}{\partial z_1} \right|
+ \cdots
+   (1 - |z_d|^2) \left| \frac{\partial \varphi}{\partial z_d} \right|
=
1 - |\varphi(z)|^2 
\end{equation*}
at every point of $\bbD^d$ if and only if it has the
form \eqref{may29a}, where $U$ is given by \eqref{may30a} with
$\cH=\cK=\bbC^d$ and, in addition, $U$ is symmetric, that is, $U$ is
equal to its transpose.

\subsection{The Schur-Agler class on the unit ball}
\label{ballcase}
Let $\Omega = \bbB_d$.  Set $\ip{z}{w} = z_1 \bar w_1 + \cdots + z_d
\bar w_d$ for any $z = (z_1,\dots,z_d)$ and $w = (w_1,\dots,w_d)$ in
$\bbC^d$, and write
\begin{equation*}
k(w,z) = \frac{1}{1 - \ip{z}{w}} , \qquad z,w\in\bbB_d.
\end{equation*}
The {\bf Arveson space} is the Hilbert space $\cH(k)$ of analytic
functions on $\bbB_d$ with reproducing kernel $k(w,z)$ (for example,
see Arveson \cite{arvesonIII} and Drury \cite{drury}).  An analytic
function $\varphi(z)$ on $\bbB_d$ is a {\bf multiplier} of the Arveson
space if multiplication by $\varphi(z)$ is an everywhere defined and
bounded operator on $\cH(k)$.  It is called a {\bf contractive
  multiplier} if this operator is a contraction, that is, it has norm
at most one.  Equivalently, $\varphi(z)$ is a contractive multiplier
of the Arveson space if and only if the kernel
$[1-\varphi(z)\overline{\varphi(w)}]/(1-\ip{z}{w})$ is nonnegative,
that is, the inequality
\begin{equation*}
  \sum_{j,k=1}^n \frac{1 - \varphi(z^{(k)})\overline{\varphi(z^{(j)})}}
                  {1 - \ip{z^{(k)}}{z^{(j)}}} c_j \bar c_k \ge 0
\end{equation*}
holds for any points $z^{(1)},\dots,z^{(n)}$ in $\bbB_d$, any complex
numbers $c_1,\dots,c_n$, and any positive integer~$n$.  For $d=1$, the
class of contractive multipliers of the Arveson space coincides with
the Schur class on $\bbB_d$.  For $d>1$, the set of contractive
multipliers of the Arveson space is a proper subset of the Schur class
(see \cite[pp.\ 99--100]{aglermccarthy} and \cite[p.\ 6]{alpaykap}).
For any $d\ge 1$, the set of contractive multipliers of the Arveson
space coincides with the class of functions which have a unitary
realization \eqref{may29a} when suitable choices are made.  In
\eqref{may30a} we choose $\cK = \cH \oplus \cdots \oplus \cH$ with $d$
summands and take
\begin{equation}\label{jun2b}
Z(z_1,\dots,z_d) = z_1 E_1 + \cdots z_d E_d,
\end{equation}
where $E_j (h_1 \oplus\cdots\oplus h_d) = h_j$, $j=1,\dots,d$.
Alternatively, write the elements of $\cK$ in column form, and let
\begin{equation*}
A = \begin{pmatrix}
  A _1 \\ \vdots \\ A_d
\end{pmatrix} ,
\qquad
B = \begin{pmatrix}
  B _1 \\ \vdots \\ B_d
\end{pmatrix} .
\end{equation*}
Then \eqref{may29a} becomes
\begin{multline*}
  \varphi(z) = D + C Z(z) (I_\cK - A Z(z))^{-1} B
  =  D + C (I_\cH - Z(z)A)^{-1} Z(z)B  \\
  =  D + C\Big(I_\cH - \sum_{j=1}^d z_j A_j \Big)^{-1} 
              \sum_{j=1}^d z_j B_j  \,.
\end{multline*}
With these choices, every function of the form \eqref{may29a} is a
contractive multiplier of the Arveson space, and every contractive
multiplier of the Arveson space has this form.  For this reason, the
class of contractive multipliers of the Arveson space is also called
the {\bf Schur-Agler class on the unit ball}.  See Ball, Trent, and
Vinnikov \cite{btv2001} and Eschmeier and Putinar \cite{ep2002}.

In \S\ref{S:calculations} of this paper, we show how to estimate the
expression \eqref{may29a} for functions on an arbitrary
domain~$\Omega$ in $\bbC^d$.  Since no extra labor is involved, the
results in \S\ref{S:calculations} are proved for operator-valued
functions.  Estimates for derivatives of complex-valued functions in
the Schur-Agler class on the polydisk and unit ball are then deduced
in \S\ref{S:polydisk} and \S\ref{S:ball} by interpreting the general
results in \S\ref{S:calculations} for these cases.  In
\S\ref{S:polydisk} and \S\ref{S:ball} we also derive special results
for the polydisk and unit ball with the aid of new Hilbert space
estimates, which appear in Lemmas \ref{dec25b} and \ref{dec27b}.
Extensions to operator-valued functions and some open questions are
discussed in \S~\ref{S:remarks}.

\section{Estimates of functions which admit a
  realization}\label{S:calculations} 
Throughout this section we assume that $\varphi(z)$ is a function of
$z = (z_1,\dots,z_d)$ on a domain $\Omega$ in $\bbC^d$ with values in
$\lfg$ for some Hilbert spaces $\cF$ and $\cG$, and that $\varphi(z)$
admits a representation in the form
\begin{equation}\label{aug23a}
\varphi(z) = D + C Z(z) (I_\cK - A Z(z))^{-1} B , \qquad z\in\Omega,
\end{equation}
where
\begin{equation}\label{aug23b}
      U=\begin{pmatrix} A & B \\[2pt]  C & D \end{pmatrix} \colon
  \begin{matrix} \cH \\[-3pt] \oplus \\[-1pt] \cF \end{matrix} \to 
  \begin{matrix} \cK \\[-3pt] \oplus \\[-1pt] \cG \end{matrix}
\end{equation}
is a unitary operator for some Hilbert spaces $\cH$ and $\cK$, and
$Z(z)$ is a function of $z = (z_1,\dots,z_d)$ in $\Omega$ with values
in $\lkh$ such that
\begin{equation}
  \label{aug23c}
  \| Z(z) \| < 1 , \qquad z\in\Omega .
\end{equation}
We also use the condition
\begin{equation}
  \label{aug23d}
 Z(z_1,\dots,z_d) = z_1 E_1 + \cdots + z_d E_d
 \quad\text{where}\quad E_j \in\lkh,\; j=1,\dots,d .
\end{equation}
The symbol $I$ always denotes an identity operator, and a subscript,
when present, indicates the underlying space.

\begin{lemma}\label{identities}
  Let $\varphi(z)$ have the form $\eqref{aug23a}$ where
  $\eqref{aug23b}$ is a unitary operator and $Z(z)$ satisfies
  $\eqref{aug23c}$.  For all $w,z \in\Omega$,
\begin{align}
  I_\cF - \varphi(z)^* & \varphi(w)  \label{may31b} \\
&=
  B^* (I_\cK -Z(z)^*A^*)^{-1} (I_\cK-Z(z)^*Z(w)) 
               (I_\cK-AZ(w))^{-1} B, \notag \\[6pt]
  I_{\cG} - \varphi(w)&\varphi(z)^*  \label{may31c}  \\
&=
  C (I_\cH - Z(w)A)^{-1} (I_\cH - Z(w)Z(z)^*) 
                   (I_\cH - A^*Z(z)^*)^{-1} C^* .  \notag
\end{align}
In particular,  $\| \varphi(z) \| \le 1$ for each $z \in\Omega$.
\end{lemma}

This result is well known, but we sketch a proof for completeness.

\begin{proof}
  To prove \eqref{may31b}, we use the identity
\begin{equation*}
  \begin{pmatrix}
    A^*A+C^*C & A^*B+C^*D \\
    B^*A+D^*C & B^*B+D^*D
  \end{pmatrix} 
=
 \begin{pmatrix}
  A^*&C^*\\ B^*&D^*
\end{pmatrix} \begin{pmatrix}
  A&B\\ C&D
\end{pmatrix}  = \begin{pmatrix}
  I&0\\0&I
\end{pmatrix} \,.
\end{equation*}
Thus
\allowdisplaybreaks{
\begin{align*}
  I - &\varphi(z)^*\varphi(w) \\
&=
  I - \Big[ D^* + B^* (I - Z(z)^* A^*)^{-1} Z(z)^* C^* \Big]
    \Big[ D + C Z(w) (I - A Z(w))^{-1} B   \Big] \\[5pt]
&=
  I-D^*D - D^*C Z(w) (I - A Z(w))^{-1} B
    - B^* (I - Z(z)^* A^*)^{-1} Z(z)^* C^*D \\*[3pt]
  &\hskip1cm
     -  B^* (I - Z(z)^* A^*)^{-1}  Z(z)^* 
                  C^*C Z(w) (I - A Z(w))^{-1} B  \\[5pt]
&=
  B^*B + B^*A Z(w) (I - A Z(w))^{-1} B
    + B^* (I - Z(z)^* A^*)^{-1} Z(z)^* A^*B \\*[3pt]
  &\hskip1cm
     -  B^* (I - Z(z)^* A^*)^{-1}  Z(z)^* 
                 (I-A^*A)  Z(w) (I - A Z(w))^{-1} B \\[5pt]
&=
  B^* (I - Z(z)^* A^*)^{-1}
  \Big\{
  (I - Z(z)^* A^*) (I - A Z(w))
  + (I - Z(z)^* A^*) AZ(w) \\*[3pt]
  &\hskip1cm
  + Z(z)^* A^* (I - A Z(w)) 
  - Z(z)^*(I-A^*A)Z(w) 
  \Big\}
  (I - A Z(w))^{-1} B \\[5pt]
&=
  B^* (I -Z(z)^*A^*)^{-1} (I-Z(z)^*Z(w)) 
               (I-AZ(w))^{-1} B ,
\end{align*}}%
proving \eqref{may31b}.  In a similar way, we have
\begin{equation*}
  \varphi(z) = D + C(I_\cH - Z(z)A)^{-1}Z(z)B 
\end{equation*}
and
\begin{equation*}
  \begin{pmatrix}
    AA^*+BB^* & AC^* + BD^* \\ 
    CA^*+DB^* & CC^*+DD^*
  \end{pmatrix} 
=
\begin{pmatrix}
  A&B\\ C&D
\end{pmatrix}  \begin{pmatrix}
  A^*&C^*\\ B^*&D^*
\end{pmatrix} = \begin{pmatrix}
  I&0\\0&I
\end{pmatrix} \,.
\end{equation*}
Parallel calculations using these identities yield \eqref{may31c}.
\end{proof}

The next lemma uses the differentiation rules
\begin{equation*}
    \frac{\partial}{\partial z_j} \big[ FG \big] 
= \frac{\partial F}{\partial z_j} G 
        + F \frac{\partial G}{\partial z_j} , \qquad
  \frac{\partial}{\partial z_j} F^{-1}
=
  - F^{-1} \frac{\partial F}{\partial z_j} F^{-1} ,
\end{equation*}
which hold for any operator-valued functions for which the operations
are defined.

\begin{lemma}
  Let $\varphi$ have the form $\eqref{aug23a}$ where
  $\eqref{aug23b}$ is a unitary operator and $Z$ satisfies
  $\eqref{aug23c}$ and $\eqref{aug23d}$.  Set
  \begin{equation*}
    L = A(I_\cH - ZA)^{-1} = (I_\cK -AZ)^{-1} A .
  \end{equation*}
  Then for each $j=1,\dots,d$,
  \begin{equation}
    \label{jun1a}
    \frac{\partial \varphi}{\partial z_j}
  =
    C (I_\cH - ZA)^{-1} E_j (I_\cK -AZ)^{-1} B .
  \end{equation}
  For every $n \ge 2$ and any $k_1,\dots,k_n$ in
  $\{1,2,\dots,d\}$,
  \begin{equation}\label{jun1b}
    \frac{\partial^{\,n} \varphi}{\partial z_{k_n} \cdots \partial
      z_{k_1}}
  =
      C (I_\cH - ZA)^{-1}
          \sum_\sigma E_{k_{\sigma(1)}} L E_{k_{\sigma(2)}} L
       \cdots  L E_{k_{\sigma(n)}} (I_\cK -AZ)^{-1} B,
  \end{equation}
  where the summation is over all permutations $\sigma$ of
  $\{1,2,\dots,n\}$.
\end{lemma}

\begin{proof}
  Note that
  \begin{align}
    \frac{\partial}{\partial z_j} (I-ZA)^{-1} 
  &=
    - (I-ZA)^{-1} (-E_jA) (I-ZA)^{-1} 
  =
    (I-ZA)^{-1} E_j L, \label{jun1c} \\
\frac{\partial}{\partial z_j} (I-AZ)^{-1}
  &=
    - (I-AZ)^{-1} (-AE_j) (I-AZ)^{-1} 
  =
     LE_j (I-AZ)^{-1}, \label{jun1d}
  \end{align}

and
\begin{align}
    \frac{\partial L}{\partial z_j}
&=
  \frac{\partial}{\partial z_j} A (I-ZA)^{-1}   \label{jun1e} \\
&\hskip1cm = 
  -A (I-ZA)^{-1} (-E_j A) (I-ZA)^{-1}
=
  L E_j L \,. \notag
\end{align}
By \eqref{jun1d},
\begin{align*}
  \frac{\partial \varphi}{\partial z_j}
&=
  CE_j (I-AZ)^{-1} B
    + CZLE_j (I-AZ)^{-1} B \\
&=
  C \big[ I+Z(I-AZ)^{-1}A  \big] E_j (I-AZ)^{-1} B \\[4pt]
&=
  C (I-ZA)^{-1} E_j (I-AZ)^{-1} B ,
\end{align*}
which is \eqref{jun1a}.  Thus by \eqref{jun1c} and \eqref{jun1d},
\begin{align*}
  \frac{\partial^{\,2} \varphi}{\partial z_{k_2} \partial z_{k_1} }
&=
  \frac{\partial }{\partial z_{k_2}}
    C (I-ZA)^{-1} E_{k_1} (I-AZ)^{-1} B \\[3pt]
&=
    C (I-ZA)^{-1} E_{k_2} L E_{k_1} (I-AZ)^{-1} B   \\[3pt]
  &\hskip1cm
    + C (I-ZA)^{-1} E_{k_1} L E_{k_2} (I-AZ)^{-1} B  \\[3pt]
&=
    C (I-ZA)^{-1}
    \; \sum_\sigma E_{k_{\sigma(1)}} L E_{k_{\sigma(2)}} \;
    (I-AZ)^{-1} B    ,
\end{align*}
where the summation is over all permutations $\sigma$ of $\{1,2\}$.
This is \eqref{jun1b} for $n=2$.

Assume that \eqref{jun1b} is true for some~$n \ge 2$.  Writing
$\sigma$ for an arbitrary permutation of $\{1,2,\dots,n\}$ and using
\eqref{jun1c}, \eqref{jun1d}, and \eqref{jun1e}, we obtain
\allowdisplaybreaks{
  \begin{align*}
    \frac{\partial^{\,n+1} \varphi}
   {\partial z_{k_{n+1}}\partial z_{k_{n}} \cdots  \partial z_{k_1} }
   \\*[0pt]
   &\hskip-2.2cm 
= 
    \frac{\partial \;\;\;} {\partial z_{k_{n+1}} }\; 
    \bigg\{  
       C (I - ZA)^{-1}
          \sum_\sigma E_{k_{\sigma(1)}} L E_{k_{\sigma(2)}} L
       \cdots  L E_{k_{\sigma(n)}} 
    (I -AZ)^{-1} B 
    \bigg\}  \\[0pt]
&\hskip-2.2cm 
= 
    C (I-ZA)^{-1} E_{k_{n+1}} L
        \sum_\sigma E_{k_{\sigma(1)}} L E_{k_{\sigma(2)}} L
       \cdots  L E_{k_{\sigma(n)}} \; (I -AZ)^{-1} B
        \\*[0pt]
  &\hskip-1.5cm
     + C (I-ZA)^{-1}
    \bigg\{ 
    \frac{\partial \;\;\;} {\partial z_{k_{n+1}} }
    \sum_\sigma E_{k_{\sigma(1)}} L E_{k_{\sigma(2)}} L
       \cdots  L E_{k_{\sigma(n)}}
       \bigg\} 
   (I -AZ)^{-1} B  \\[0pt]
  &\hskip-1.5cm
     + C (I-ZA)^{-1} \; 
     \sum_\sigma E_{k_{\sigma(1)}} L E_{k_{\sigma(2)}} L
       \cdots  L E_{k_{\sigma(n)}}
    L E_{k_{n+1}} (I -AZ)^{-1} B   \\[0pt]
  &\hskip-2.2cm 
= 
  C (I-ZA)^{-1} 
   \bigg\{
   \Big[E_{k_{n+1}} L\Big] 
       \sum_\sigma E_{k_{\sigma(1)}} L E_{k_{\sigma(2)}} L
       \cdots  L E_{k_{\sigma(n)}}    \\*[0pt]
  &\hskip-.9cm 
   + \sum_\sigma E_{k_{\sigma(1)}} \Big[LE_{k_{n+1}}L\Big]  
                 E_{k_{\sigma(2)}} L
       \cdots  L E_{k_{\sigma(n)}}     \\*[0pt]
  &\hskip-.9cm 
   + \sum_\sigma E_{k_{\sigma(1)}} L E_{k_{\sigma(2)}} 
                \Big[LE_{k_{n+1}}L\Big] 
       \cdots  L E_{k_{\sigma(n)}}      \\*[0pt]
  &\hskip-.9cm  
   + \cdots   
   + \sum_\sigma E_{k_{\sigma(1)}} L E_{k_{\sigma(2)}} L
       \cdots  \Big[LE_{k_{n+1}}L\Big]  E_{k_{\sigma(n)}}      
                \\[0pt]
  &\hskip.5cm 
   +  \sum_\sigma E_{k_{\sigma(1)}} L E_{k_{\sigma(2)}} L
       \cdots  L E_{k_{\sigma(n)}} \Big[L E_{k_{n+1}}\Big]
   \bigg\}  (I-AZ)^{-1} B .
  \end{align*}}%
  This has the form \eqref{jun1b} with $n$ replaced by $n+1$.
  By induction, \eqref{jun1b} holds for all~$n \ge 2$.
\end{proof}

The norm estimates in Lemmas \ref{dec12c} and \ref{dec12d} are
understood to be in the pointwise sense.

\begin{lemma}\label{dec12c}
  Let $\varphi$ have the form $\eqref{aug23a}$ where $\eqref{aug23b}$
  is a unitary operator and $Z$ satisfies $\eqref{aug23c}$ and
  \eqref{aug23d}.  For each $j=1,\dots,d$,
  \begin{align}
    \| E_j (I_\cK - AZ)^{-1} B \|
&\le
  \| I_\cF - \varphi^*\varphi   \|^{\frac12}
    \| E_j (I_\cK - Z^*Z)^{-1}E_j^* \|^{\frac12} ,
    \label{dec12e} \\[3pt]
    \| C(I_\cH - ZA)^{-1}E_j \|
&\le
   \| I_\cG - \varphi\varphi^* \|^{\frac12}
    \| E_j^* (I_\cH - ZZ^*)^{-1}E_j \|^{\frac12} .
  \label{dec12f}
  \end{align}
Moreover,
  \begin{align}
    \| (I_\cK - AZ)^{-1} B \|
&\le
  \sqrt{\frac{\| I_\cF - \varphi^*\varphi   \|}
          {1 - \| Z \|^2}} \,,     \label{jun1f} \\[8pt]
    \| C(I_\cH - ZA)^{-1} \|
&\le
  \sqrt{\frac{\| I_{\cG} - \varphi\varphi^*   \|}
          {1 - \| Z \|^2}} \,.    \label{jun1g}
  \end{align}
\end{lemma}

\begin{proof}
By \eqref{may31b}, $\| \big( I - \varphi^*\varphi
  \big)^\frac12 \| = \| \big( I-Z^*Z\big)^\frac12
  (I-AZ)^{-1}B\|$.  Therefore
  \begin{align*}
    \|E_j(I - A Z)^{-1} B \| 
 &=
  \|E_j \big( I-Z^*Z\big)^{-\frac12}
      \big( I-Z^*Z\big)^\frac12 (I-AZ)^{-1}B  \| \\[4pt]
&\le
    \|E_j \big( I-Z^*Z\big)^{-\frac12}  \|\;
  \| I - \varphi^*\varphi \|^\frac12 \\[4pt]
&=
  \| I - \varphi^*\varphi   \|^{\frac12}
    \| E_j (I - Z^*Z)^{-1}E_j^* \|^{\frac12} \,,
  \end{align*}
  which is \eqref{dec12e}.  Similarly, using \eqref{may31c} we obtain
  \allowdisplaybreaks{
  \begin{align*}
    \| C(I -  Z  A)^{-1} E_j\|
&= \|E_j^*(I-A^*Z^*)^{-1} C^* \| \\*[6pt]
&=
  \| E_j^* \big( I-ZZ^*\big)^{-\frac12}  \big(
  I-ZZ^*\big)^\frac12
   (I-A^*Z^*)^{-1} C^* \|  \\[6pt]
&\le
   \| E_j^* \big( I-ZZ^*\big)^{-\frac12}\|\;
   \| \big( I - \varphi \varphi^* \big)^\frac12 \|  \\[6pt]
&=
   \| I - \varphi\varphi^*   \|^{\frac12}
    \| E_j^* (I - ZZ^*)^{-1}E_j \|^{\frac12} \,,
  \end{align*}}%
which is \eqref{dec12f}.  To prove \eqref{jun1f} and \eqref{jun1g} we
proceed in the same way without the factors $E_j$ and $E_j^*$, and at
the last stage we use the inequality
  \begin{equation*}
    \| \big( I-Z^*Z\big)^{-\frac12}  \|^2
= \| \big( I-Z^*Z\big)^{-1} \|
\le \sum_{k=0}^\infty \| Z^*Z \|^k = \frac{1}{1 - \| Z \|^2} 
  \end{equation*}
  and its companion $\| \big( I-ZZ^*\big)^{-\frac12} \|^2 \le
  1/(1 - \| Z \|^2)$.
\end{proof}

\goodbreak
\begin{lemma}\label{dec12d}
  Let $\varphi$ have the form $\eqref{aug23a}$ where
  $\eqref{aug23b}$ is a unitary operator and $Z$ satisfies
  $\eqref{aug23c}$ and $\eqref{aug23d}$.

\noindent $(i)$  For each  $j=1,\dots,d$,
    \begin{align*}
      \left\|
        \frac{\partial\varphi}{\partial z_j}
      \right\|
&\le
  \frac{\|  I_\cF - \varphi^*\varphi \|^\frac12 
        \| I_\cG - \varphi\varphi^*\|^\frac12}
      {(1 - \|Z\|^2)^\frac12} \cdot \\
&\hskip2cm
    \cdot \min
    \left\{
       \| E_j (I_\cK - Z^*Z)^{-1}E_j^* \|^{\frac12} ,
       \| E_j^* (I_\cH - ZZ^*)^{-1}E_j \|^{\frac12} 
    \right\}  \,.
    \end{align*}

\noindent $(ii)$ 
   For every $n \ge 2$ and any $k_1,\dots,k_n$ in
    $\{1,2,\dots,d\}$,
    \begin{align*}
       \left\|
       \frac{\partial^{\,n} \varphi}{\partial z_{k_n} \cdots \partial
      z_{k_1}}
    \right\|
&\le
    (n-2)! \;
    \frac{\|  I_\cF - \varphi^*\varphi \|^\frac12 
        \| I_\cG - \varphi\varphi^*\|^\frac12}
      {(1 - \|Z\|)^{n-1}}\; \cdot \\
&\hskip1cm
    \cdot 
    \sum_{\substack{p,q=1 \\ p \neq q}}^n
   \| E_{k_p} (I_\cK - Z^*Z)^{-1}E_{k_p}^* \|^{\frac12} \,
       \| E_{k_q}^* (I_\cH - ZZ^*)^{-1}E_{k_q} \|^{\frac12}  \;.
    \end{align*}
\end{lemma}

\begin{proof}
  (i) By \eqref{jun1a},  \eqref{dec12e}, and \eqref{jun1g},
  \begin{align*}
     \left\|
        \frac{\partial\varphi}{\partial z_j}
      \right\|
&\le
  \| C (I - ZA)^{-1} \|\, \| E_j (I -AZ)^{-1} B \| \\
&\le
\frac{\|  I - \varphi^*\varphi \|^\frac12 
        \| I - \varphi\varphi^*\|^\frac12}
      {(1 - \|Z\|^2)^\frac12}
      \| E_j (I - Z^*Z)^{-1}E_j^* \|^{\frac12} .
  \end{align*}
Using \eqref{dec12f} and \eqref{jun1f}, we get
  \begin{align*}
     \left\|
        \frac{\partial\varphi}{\partial z_j}
      \right\|
&\le
  \| C (I - ZA)^{-1} E_j \|\, \| (I -AZ)^{-1} B \| \\
&\le
\frac{\|  I - \varphi^*\varphi \|^\frac12 
        \| I - \varphi\varphi^*\|^\frac12}
      {(1 - \|Z\|^2)^\frac12}
      \| E_j^* (I - ZZ^*)^{-1}E_j \|^{\frac12} .
  \end{align*}
  The inequality in (i) follows.

(ii) By \eqref{jun1b}, \eqref{dec12e}, and \eqref{dec12f},
\begin{align*}
  \left\|
    \frac{\partial^{\,n} \varphi}{\partial z_{k_n} \cdots \partial
      z_{k_1}}
  \right\| \\[4pt]
  &\hskip-2cm =
       \big\|
      \sum_\sigma C (I - ZA)^{-1}E_{k_{\sigma(1)}}
L E_{k_{\sigma(2)}} L
       \cdots E_{k_{\sigma(n-1)}} L  
E_{k_{\sigma(n)}}(I -AZ)^{-1} B   \big\| \\[6pt]
&\hskip-2cm \le
  \sum_\sigma\big\| C (I - ZA)^{-1} 
             E_{k_{\sigma(1)}} \big\| \, \|L\|^{n-1}
   \| E_{k_{\sigma(n)}}(I -AZ)^{-1} B \| \\[6pt]
&\hskip-2cm \le
  \|L\|^{n-1} \sum_\sigma
 \| I - \varphi\varphi^* \|^{\frac12}
    \| E_{{k_{\sigma(1)}}}^* (I - ZZ^*)^{-1}
               E_{{k_{\sigma(1)}}} \|^{\frac12}
   \\
&\hskip2cm
   \cdot
 \| I - \varphi^*\varphi   \|^{\frac12}
    \| E_{k_{\sigma(n)}} (I - Z^*Z)^{-1}
               E_{k_{\sigma(n)}}^* \|^{\frac12} \,.
\end{align*}
Here the summations run over all permutations $\sigma$ of
$\{1,\dots,n\}$.  Since
  \begin{equation}\label{dec26c}
    \| L \| = \Big\| \sum_{k=0}^\infty A [ZA]^k \Big\|
    \le \sum_{k=0}^\infty \|Z\|^k 
    = \frac{1}{1-\|Z\|} \;,
  \end{equation}
we obtain
\allowdisplaybreaks{
\begin{align*}
    \left\|
    \frac{\partial^{\,n} \varphi}{\partial z_{k_n} \cdots \partial
      z_{k_1}}
  \right\|
&\le
   \frac{\|  I - \varphi^*\varphi \|^\frac12
        \| I - \varphi\varphi^*\|^\frac12}
      {(1 - \|Z\|)^{n-1}}\; \cdot \\*
&\hskip.75cm
   \cdot \sum_\sigma
    \| E_{{k_{\sigma(1)}}}^* (I - ZZ^*)^{-1}
               E_{{k_{\sigma(1)}}} \|^{\frac12}
     \| E_{k_{\sigma(n)}} (I - Z^*Z)^{-1}
               E_{k_{\sigma(n)}}^* \|^{\frac12} \\
&=
  (n-2)! \;
    \frac{\|  I - \varphi^*\varphi \|^\frac12
        \| I - \varphi\varphi^*\|^\frac12}
      {(1 - \|Z\|)^{n-1}}\; \cdot \\*
&\hskip.75cm
    \cdot
    \sum_{\substack{p,q=1 \\ p \neq q}}^n
   \| E_{k_p} (I - Z^*Z)^{-1}E_{k_p}^* \|^{\frac12} \,
       \| E_{k_q}^* (I - ZZ^*)^{-1}E_{k_q} \|^{\frac12}  \;.
\end{align*}}%
This proves the inequality in (ii).
\end{proof}

\section{The Schur-Agler class on the polydisk}
\label{S:polydisk}

To begin, we apply the esimates derived in \S\ref{S:calculations} to
complex-valued functions on the polydisk~$\bbD^d$.  Set
\begin{equation*}
  {\|}z{\|}_\infty = \max_{j=1,\dots,d}  |z_j| \;,
\qquad z = (z_1,\dots,z_d) \in \bbC^d \,.
\end{equation*}

\begin{theorem}\label{dec13a}
  Let $\varphi(z)$ be a complex-valued function in the Schur-Agler
  class on the polydisk~$\bbD^d$.  For each $z\in\bbD^d$ and
  $j=1,\dots,d$, 
  \begin{equation}
    \label{dec13b}
    \left|
      \frac{\partial \varphi}{\partial z_j}
    \right|
\le
    \frac{1 - |\varphi(z)|^2}
    {\sqrt{1-|z_j|^2}\sqrt{1-\|z\|_\infty^2}\vertspace}  \;.
  \end{equation}
  For each $z\in\bbD^d$, every $n \ge 2$, and any $k_1,\dots,k_n$ in
  $\{1,2,\dots,d\}$,
    \begin{equation}\label{dec13c}
       \left|
       \frac{\partial^{\,n} \varphi}{\partial z_{k_n} \cdots \partial
      z_{k_1}}
    \right|
\le
   (n-2)! \; \frac{1-|\varphi(z)|^2}
              {(1-\|z\|_\infty)^{n-1}}
    \sum_{\substack{p,q=1 \\ p \neq q}}^n
    \frac{1}{\sqrt{1-|z_{k_p}|^2}\sqrt{1-|z_{k_q}|^2}} \;.
    \end{equation}
\end{theorem}

 \begin{proof}
   Represent $\varphi(z)$ in the form \eqref{aug23a} with $\cF = \cG =
   \bbC$, the complex numbers in the Euclidean metric, and choose
   $\Omega = \bbD^d$.  Following \S\ref{polydiskcase}, we assume that
   the Hilbert spaces $\cH$ and $\cK$ satisfy $\cK = \cH = \cH_1
   \oplus \cdots \oplus \cH_d$\,, and that $Z(z) = z_1 E_1 + \cdots +
   z_d E_d$, where $E_j$ is the orthogonal projection onto $\cH_j$,
   $j=1,\dots,d$.  Then $\|Z(z)\| \le \|z\|_\infty$.  For any
   $j=1,\dots,d$,
   \begin{equation*}
    E_j (I - Z(z)Z(z)^*)^{-1} E_j^* = E_j^* (I - Z(z)^*Z(z))^{-1} E_j
  =
    \frac{1}{1 - |z_j|^2} \; E_j \,,
   \end{equation*}
   and therefore
 \begin{equation*}
   \| E_j (I - Z(z)Z(z)^*)^{-1} E_j^* \|^\frac12
= \| E_j^* (I - Z(z)^*Z(z))^{-1} E_j \|^\frac12
\le
   \frac{1}{\sqrt{1 - |z_j|^2}} \,.
 \end{equation*}
 Now apply Lemma~\ref{dec12d}.
\end{proof}

Recall that the Schur-Agler class on $\bbD^2$ is the full Schur class
of analytic functions which are defined and bounded by one on
$\bbD^2$.

\begin{corollary}
  Let $\varphi(z)$ be any complex-valued function which is analytic
  and bounded by one on~$\bbD^2$. For each $z\in\bbD^2$ and $j=1,2$,
  \begin{equation}
    \label{dec16b}
    \left|
      \frac{\partial \varphi}{\partial z_j}
    \right|
\le
    \frac{1 - |\varphi(z)|^2}
    {\sqrt{1-|z_j|^2}\sqrt{1-\|z\|_\infty^2}\vertspace}  \;.
  \end{equation}
  For each $z\in\bbD^2$ and any nonnegative integers $n_1,n_2$
  with $n = n_1 + n_2 \ge 2$,
  \begin{multline}\label{dec16a}
    \left|
      \frac{\partial^{\,n} \varphi}
    {\partial z_2^{n_2} \partial z_1^{n_1}}
    \right|
\le
   (n-2)! \; \frac{1-|\varphi(z)|^2}
              {(1-\|z\|_\infty)^{n-1}} \cdot \\
   \cdot \bigg[ \frac{n_1^2-n_1}{1-|z_1|^2}
       + \frac{2n_1 n_2}{\sqrt{1-|z_1|^2}\sqrt{1-|z_2|^2}}
       +  \frac{n_2^2-n_2}{1-|z_2|^2} \bigg] \;.
  \end{multline}
\end{corollary}

\begin{proof}
  The inequality \eqref{dec16b} is a repeat of \eqref{dec13b}.  To
  deduce \eqref{dec16a} from \eqref{dec13c}, write
  \begin{equation*}
   (k_1,\dots,k_n) = 
    (\, \overbrace{1,\dots,1}^{n_1},
       \overbrace{2,\dots,2}^{n_2} ) \,,
  \end{equation*}
  and collect similar terms.
\end{proof}

\goodbreak
The next result takes advantage of the special form of the polydisk.

\begin{theorem}\label{dec20a}
  Let $\varphi(z)$ be a complex-valued function in the Schur-Agler
  class on the polydisk $\bbD^d$.  Then for each $z\in\bbD^d$ and any
  nonnegative integers $n_1,\dots,n_d$ which are not all zero,
  \begin{equation}\label{dec20b}
    \left|
      \frac{\partial^{\,n} \varphi}
    {\partial z_d^{n_d} \cdots \partial z_1^{n_1}}
    \right|
\le
    n_1! \cdots n_d!\; \frac{1-|\varphi(z)|^2}
     {(1-\|z\|_\infty^2)(1-\|z\|_\infty)^{n-1}} \,,
  \end{equation}
  where $n = n_1 + \cdots + n_d$.
\end{theorem}

The weaker inequality 
  \begin{equation}\label{dec20c}
    \left|
      \frac{\partial^{\,n} \varphi}
    {\partial z_d^{n_d} \cdots \partial z_1^{n_1}}
    \right|
\le
    n! \; \frac{1-|\varphi(z)|^2}
     {(1-\|z\|_\infty^2)(1-\|z\|_\infty)^{n-1}} \,,
  \end{equation}
  follows from Theorem~\ref{dec13a}.  In fact, when $n=1$ this is
  immediate from \eqref{dec13b}.  Suppose $n \ge 2$.  Set
  \begin{equation}\label{dec17a}
   (k_1,\dots,k_n) = 
    (\, \overbrace{1,\dots,1}^{n_1},
       \overbrace{2,\dots,2}^{n_2},
       \dots,
       \overbrace{d,\dots,d}^{n_d}\, ) \,.
  \vphantom{\underbrace{x}}
  \end{equation}
  By \eqref{dec13c},
  \begin{multline*}
    \left|
      \frac{\partial^{\,n} \varphi}
    {\partial z_d^{n_d} \cdots \partial z_1^{n_1}}
    \right|
=
  \left|
       \frac{\partial^{\,n} \varphi}{\partial z_{k_n} \cdots \partial
      z_{k_1}}
    \right|
\le
  (n-2)! \; \frac{1-|\varphi(z)|^2}
              {(1-\|z\|_\infty)^{n-1}}
    \sum_{\substack{p,q=1 \\ p \neq q}}^n
    \frac{1}{1-\|z\|_\infty^2} \\
=
  (n-2)! \; \frac{1-|\varphi(z)|^2}
     {(1-\|z\|_\infty^2)(1-\|z\|_\infty)^{n-1}} \; (n^2-n)
=
     n! \; \frac{1-|\varphi(z)|^2}
          {(1-\|z\|_\infty^2)(1-\|z\|_\infty)^{n-1}} \,,
  \end{multline*}
  yielding \eqref{dec20c}.  
  
  A variation on the method, which is valid for the polydisk case but
  not generally, is needed to prove \eqref{dec20b}.

  \begin{lemma}\label{dec25b}
    Let $L\in\lh$, where $\cH$ is a Hilbert space.  Assume that $\cH =
    \cH_1 \oplus\cdots\oplus \cH_d$\,, and let $E_1,\dots,E_d$ be the
    orthogonal projections onto the summands.  Suppose that $n = n_1 +
    \cdots + n_d$ for some nonnegative integers $n_1,\dots,n_d$.
    Assume that $n \ge 2$, and set
    \begin{equation}\label{dec26a}
      K = \sum E_{j_1} L E_{j_2} L \cdots E_{j_{n-1}} L E_{j_n},
    \end{equation}
    where the summation is over all distinct arrangements
    $(j_1,j_2,\dots,j_n)$ of the tuple
    \begin{equation}\label{dec26b}
      (\,
    \overbrace{1,\dots,1}^{n_1}, \overbrace{2,\dots,2}^{n_2}, \dots,
    \overbrace{d,\dots,d}^{n_d}\, ) \;. 
    \vphantom{\underbrace{x}{x}}
    \end{equation}
    Then
    \begin{equation}
      \label{dec25c}
      \| K \| \le \| L \|^{n-1} \, .    \end{equation}
  \end{lemma}
  
  The number of terms in \eqref{dec26a} is the multinomial coefficient
   \begin{equation*}
     (n;n_1,n_2,\dots,n_d) = \frac{n!}{n_1! n_2! \cdots n_d!} \,.
   \end{equation*}
   Lemma~\ref{dec25b} is not stated in maximum generality, but it will
   suit our application.

  \begin{proof}[Proof of Lemma~\ref{dec25b}]
    We first illustrate the method for the example $n_1 = 3$, $n_2 =
    2$, and $d=2$.  Write
    \begin{equation}\label{dec25a}
      K = K[1,1,1,2,2] = \sum
                E_{j_1} L E_{j_2} L  E_{j_{3}} L E_{j_4} L E_{j_5},
    \end{equation}
    where the summation is over all distinct arrangements of
    $(1,1,1,2,2)$:
    \begin{align*}
      &(1,1,1,2,2), (1,1,2,1,2), (1,1,2,2,1), (1,2,1,1,2),
    (1,2,1,2,1), (1,2,2,1,1) , \\
    &(2,1,1,1,2), (2,1,1,2,1), (2,1,2,1,1), (2,2,1,1,1)\;.
    \end{align*}
    Define $K[1,1,2,2], K[1,2,2], \dots$ in a similar way.  For
    example, $K[1,1,2,2]$ is the sum of all terms $E_{j_1} L E_{j_2} L
    E_{j_{3}} L E_{j_4}$ over the set of distinct arrangements of
    $(1,1,2,2)$.  Now split the sum \eqref{dec25a} into two parts, the
    first part consisting of all terms with $j_1 = 1$ and the second
    consisting of all terms with $j_1 = 2$.  The result is
    \begin{equation*}
      K =  E_1 L K[1,1,2,2] + E_2 L K[1,1,1,2].
    \end{equation*}
    Since $E_1$ and $E_2$ have orthogonal ranges, for any $f \in\cH$,
    \begin{align*}
      \| Kf \|^2 &= 
    \| E_1 L K[1,1,2,2]f \|^2 + \| E_2 L K[1,1,1,2] f \|^2 \\
&\le
    \|L \|^2 \Big\{ \|K[1,1,2,2]f \|^2 + \| K[1,1,1,2] f\|^2 \Big\} \,.
    \end{align*}
    Next decompose $K[1,1,2,2]f$ and $K[1,1,1,2] f$ in the same way:
    \begin{align*}
       \| Kf \|^2 
&\le 
      \|L \|^2 \Big\{
       \| E_1LK[1,2,2]f\|^2 + \| E_2LK[1,1,2]f\|^2 \\
&\hskip2cm
       + \|E_1LK[1,1,2]f \|^2 + \| E_2LK[1,1,1]f\|^2
        \Big\} .
    \end{align*}
    The two middle terms in the last expression combine and
    give
    \begin{align*}
      \| E_2LK[1,1,2]f\|^2 +   \|E_1LK[1,1,2]f \|^2
&=
    \| (E_2+E_1)LK[1,1,2]f\|^2 \\
&\le 
    \| L \|^2 \| K[1,1,2]f \|^2 ,
    \end{align*}
   because $E_1 + E_2$ is also an orthogonal projection.  Thus
   \begin{equation*}
  \| Kf \|^2 
\le
   \|L \|^4 \Big\{ \|K[1,2,2]f \|^2 + \|K[1,1,2]f \|^2 
       + \|K[1,1,1]f \|^2 \Big\} .
   \end{equation*}
   Continuing in this way, we obtain
\allowdisplaybreaks{
   \begin{align*}
     \| Kf \|^2 
&\le 
   \|L \|^4  \Big\{
  \| E_1 L K[2,2]f \|^2 +  \| E_2 L K[1,2]f \|^2 \\*[3pt]
&\hskip2cm
  + \| E_1 L K[1,2]f \|^2 +  \| E_2 L K[1,1]f \|^2 
  + \| E_1 L K[1,1]f \|^2
    \Big\} \\[5pt]
&\le
  \|L \|^6 \Big\{ 
  \| K[2,2]f \|^2 + \| K[1,2]f \|^2 + \| K[1,1]f \|^2
  \Big\} \\[5pt]
&=
   \|L \|^6 \Big\{ 
   \| E_2 LE_2 f \|^2 
    + \Big( \| E_1 LE_2 f \|^2 + \| E_ 2LE_1 f \|^2 \Big)
    + \| E_1 LE_1 f \|^2
   \Big\} \\[5pt]
&\le
    \|L \|^8 \Big\{ 
   \| E_2 f \|^2 + \| E_1 f \|^2
  \Big\} \\[5pt]
&=
   \|L \|^8\, \| f \|^2 .
   \end{align*}}%
   Thus $\| K \| \le \|L \|^4$, which proves the lemma for the
   example.
   
   In the general case, we consider any $f\in\cH$ and write
   \begin{align*}
Kf
&=
     E_1 L K[\,\overbrace{1,\dots,1}^{n_1-1},
     \overbrace{2,\dots,2}^{n_2}, 
       \dots, \overbrace{d,\dots,d}^{n_d}\,]f \\
&\hskip3cm
    + E_2 L K[\,\overbrace{1,\dots,1}^{n_1}, 
        \overbrace{2,\dots,2}^{n_2-1},
     \dots, \overbrace{d,\dots,d}^{n_d}\,]f  \\
&\hskip3cm
    + \cdots
    + E_d L K[\,\overbrace{1,\dots,1}^{n_1}, 
        \overbrace{2,\dots,2}^{n_2},
     \dots, \overbrace{d,\dots,d}^{n_d-1}\,]f \;.
   \end{align*}
   As in the example, we use orthogonality and repeatedly break down
   all expressions into smaller and smaller patterns which have the
   same length at each stage.  Whenever a pattern
   $K[p_1,\dots,p_\ell]$ appears more than once, it is in a group
   \begin{equation*}
     E_{q_1}LK[p_1,\dots,p_\ell]f ,\; E_{q_2}LK[p_1,\dots,p_\ell]f ,
     \, \dots\,, E_{q_m}LK[p_1,\dots,p_\ell]f 
   \end{equation*}
   with distinct $q_1,q_2,\dots,q_m$ in $\{1,\dots,d\}$.  Since these
   expressions are orthogonal, the group can be collapsed.  In the
   end, we obtain $\| Kf \|^2 \le \| L \|^{2(n-1)}\, \| f \|^2$.  By
   the arbitrariness of $f$, \eqref{dec25c} follows.
  \end{proof}

  \begin{proof}[Proof of Theorem~\ref{dec20a}]
    We use the representation \eqref{aug23a} of $\varphi(z)$ with the
    same choices of Hilbert spaces and operators as in the proof of
    Theorem~\ref{dec13a}.  When $n=1$, \eqref{dec20b} is immediate
    from \eqref{jun1a}, \eqref{jun1f}, and \eqref{jun1g}.  Suppose
    $n\ge 2$.  By \eqref{jun1b} applied with $(k_1,k_2,\dots,k_n)$ as
    in \eqref{dec17a},
    \begin{align}
    &\frac{\partial^{\,n} \varphi}
    {\partial z_d^{n_d} \cdots \partial z_1^{n_1}} \label{dec27a} \\
&\hskip.7cm =
      \frac{\partial^{\,n} \varphi}{\partial z_{k_n} \cdots \partial
      z_{k_1}}  \notag \\[4pt]
  &\hskip.7cm =
      C (I_\cH - Z(z)A)^{-1} \cdot
          \sum_\sigma E_{k_{\sigma(1)}} L(z) E_{k_{\sigma(2)}} L(z)
       \cdots  L(z) E_{k_{\sigma(n)}} \cdot \notag \\
  &\hskip5cm
    \cdot (I_\cK -AZ(z))^{-1} B \notag \\[8pt]
&\hskip.7cm =
   n_1! n_2! \cdots n_d! \;
   C (I_\cH - Z(z)A)^{-1} \cdot
          \sum E_{j_1} L(z) E_{j_2} L(z)
       \cdots  L(z) E_{j_n} \cdot \notag \\[5pt]
  &\hskip5cm
    \cdot (I_\cK -AZ(z))^{-1} B \;,   \notag
    \end{align}
    where the last summation is over all distinct arrangements
    $(j_1,j_2,\dots,j_n)$ of the tuple \eqref{dec26b}.  The factorials
    $n_1!$, $n_2!$, \dots, $n_d!$ account for the permutations of the
    groups in \eqref{dec17a}.  By Lemma~\ref{dec25b} applied with $L =
    L(z)$,
    \begin{equation*}
      \left|
     \frac{\partial^{\,n} \varphi}
    {\partial z_d^{n_d} \cdots \partial z_1^{n_1}}
      \right|
\le
   \big\| C (I_\cH - Z(z)A)^{-1}\big\|  \big\| L(z) \big\|^{n-1}  
          \big\|(I_\cK -AZ(z))^{-1} B\big\| \,.
    \end{equation*}
    Then \eqref{dec20b} follows from \eqref{jun1f}, \eqref{jun1g}, and
    \eqref{dec26c}.
  \end{proof}
  
\goodbreak
  A theorem of F.\ Wiener \cite{BDK, Bohr, landaugaier} asserts that
  if $\varphi(z)$ is a complex-valued function which is analytic and
  bounded by one on $\bbD$, and if $\varphi(z) = \sum_{k=0}^\infty c_k
  z^k$, then $|c_k| \le 1 - |c_0|^2$ for every $k\ge 1$.  An immediate
  consequence of Theorem~\ref{dec20a} is an analog for the Schur-Agler
  class on the polydisk.

  \begin{corollary}
    Let $\varphi(z)$ be a complex-valued function in the Schur-Agler
    class on the polydisk~$\bbD^d$.  If
    \begin{equation*}
      \varphi(z) = \sum_{k_1,\dots,k_d=0}^\infty c_{k_1\cdots k_d}
           z_1^{k_1} \cdots z_d^{k_d},
    \end{equation*}
   then $|c_{n_1\cdots n_d}| \le 1 - |c_{0\cdots 0}|^2$
   for all nonnegative integers $n_1,\dots,n_d$ which are not all
   zero.
  \end{corollary}
  
  The question arises if one of the estimates \eqref{dec13c} or
  \eqref{dec20b} is always better than the other.  This is not the
  case.  Take $d=2$, $n \ge 2$, and use the equivalent form
  \eqref{dec16a} for \eqref{dec13c}.  For $z=(0,0)$, \eqref{dec20b} is
  a better estimate than \eqref{dec16a} except in trivial cases.
  However, it is not hard to see that when $|z_2| > |z_1|$ and
  $n_2=0$, the reverse is true.

\section{The Schur-Agler class on the 
unit ball}\label{S:ball}

In a similar way we can specialize the results of
\S\ref{S:calculations} to the Schur-Agler class on the unit ball
$\bbB_d$.  For any $z = (z_1,\dots,z_d)$ in $\bbC^d$, write
\begin{equation*}
  \| z \|_2 = \ipzz^\frac12 = \sqrt{|z_1|^2 + \cdots + |z_d|^2}\; 
\end{equation*}
and
\begin{equation*}
  \hat z_j = (z_1,\dots,z_{j-1},0,z_{j+1},\dots,z_d) .
\end{equation*}

\begin{theorem}\label{dec13d}
  Let $\varphi(z)$ be a complex-valued function in the Schur-Agler
  class on the unit ball.  For each $z\in \bbB_d$ and any
  nonnegative integers $n_1,\dots,n_d$ which are not all zero,
    \begin{equation}
      \label{dec13f}
             \left|
       \frac{\partial^{\,n} \varphi}{\partial z_{d}^{n_d} \cdots \partial
      z_1^{n_1}}
    \right|
\le
  (n-1)!\, \frac{1-|\varphi(z)|^2}
          {(1- \|z\|_2^2 )(1- \|z\|_2)^{n-1}} \,
     \sum_{j=1}^d n_j \sqrt{1-\|\hat z_j\|_2^2} \;,
    \end{equation}
   where $n=n_1 + \cdots + n_d$.
\end{theorem}

\begin{proof}
  Represent $\varphi(z)$ in the form \eqref{aug23a} with $\cF = \cG =
  \bbC$, and let $\Omega = \bbB_d$.  As in \S\ref{ballcase}, we take
  $\cK = \cH \oplus \cdots \oplus \cH$ with $d$ summands, and we write
  the elements of $\cK$ in column form.  Then
  \begin{equation*}
    Z(z) =
  \begin{bmatrix} z_1 I_\cH & \cdots & z_d I_\cH \end{bmatrix}
     \qquad\text{and}\qquad
     E_j = \begin{bmatrix} 0 &\cdots &0 &I_\cH &0 &\cdots &0
            \end{bmatrix},
  \end{equation*}
  where the identity operator appears in the $j$-th position in the
  formula for $E_j$, $j=1,\dots,d$.  Thus $Z(z)Z(z)^* = \|z\|_2^2 \,
  I_\cH$, and hence $\| Z(z) \| = \| Z(z)Z(z)^* \|^\frac12 = \|z\|_2$
  and
  \begin{equation}\label{dec17b}
    \| E_j^* (I-Z(z)Z(z)^*)^{-1}E_j \| = \frac{1}{1-\|z\|_2^2} \,.
  \end{equation}
  Using $E_jZ(z)^* = \bar z_j I_\cH$ and $E_jE_j^* =
  I_\cH$, we get
  \begin{align*}
    E_j (I_\cK-Z(z)^*Z(z))^{-1} E_j^* \\ 
&\hskip-2.5cm =
  E_j \big(
   I_\cK + Z(z)^*Z(z) + Z(z)^*Z(z)Z(z)^*Z(z) + \cdots
   \big) E_j^* \\[5pt]
&\hskip-2.5cm =
  I_\cH + \bar z_j I_\cH z_j
   + \bar z_j Z(z)Z(z)^* z_j
   + \bar z_j Z(z)Z(z)^*Z(z)Z(z)^* z_j + \cdots \\[5pt]
&\hskip-2.5cm =
  I_\cH + |z_j|^2 (I_\cH - Z(z)Z(z)^*)^{-1} \\
&\hskip-2.5cm =
\left(
  1 + \frac{|z_j|^2}{1 - \|z\|_2^2}
\right) I_\cH \;.  
  \end{align*}
Therefore 
\begin{equation}\label{dec17c}
  \| E_j (I_\cK-Z(z)^*Z(z))^{-1} E_j^* \|
=
  1 + \frac{|z_j|^2}{1 - \|z\|_2^2}  
=
  \frac{1 - \| \hat z_j \|_2^2}{1 - \|z\|_2^2} \;.
\end{equation}
By Lemma \ref{dec12d}(i), 
\begin{equation*}
  \left|
    \frac{\partial\varphi}{\partial z_j}
  \right|
\le
  \frac{1-|\varphi(z)|^2}{\sqrt{1-\|z\|_2^2}} 
   \; \sqrt{\frac{1 - \| \hat z_j \|_2^2}{1 - \|z\|_2^2}}
=
   \frac{1-|\varphi(z)|^2}{1-\|z\|_2^2} 
   \sqrt{1 - \| \hat z_j \|_2^2} \;,
\end{equation*}
which is \eqref{dec13f} when $n=1$.

Consider $n \ge 2$ and any $k_1,\dots,k_n$ in $\{1,2,\dots,d\}$.  By
Lemma \ref{dec12d}(ii), \eqref{dec17b}, and \eqref{dec17c},
\allowdisplaybreaks{
\begin{align*}
   \left|
       \frac{\partial^{\,n} \varphi}{\partial z_{k_n} \cdots \partial
      z_{k_1}}
    \right|
&\le 
  (n-2)! \; \frac{1-|\varphi(z)|^2}{(1-\|z\|_2)^{n-1}} 
  \sum_{\substack{p,q=1 \\ p \neq q}}^n
  \frac{1}{\sqrt{1-\|z\|_2^2}}
  \sqrt{\frac{1-\|\hat z_{k_q}\|_2^2}{1-\|z\|_2^2}} \\
&=
  (n-2)! \; \frac{1-|\varphi(z)|^2}
   {(1-\|z\|_2^2)(1-\|z\|_2)^{n-1}} 
  \sum_{\substack{p,q=1 \\ p \neq q}}^n
  \sqrt{1-\|\hat z_{k_q}\|_2^2} \\
&=
  (n-1)!\; \frac{1-|\varphi(z)|^2}
          {(1- \|z\|_2^2 )(1- \|z\|_2)^{n-1}} \,
     \sum_{j=1}^n \sqrt{1-\|\hat z_{k_j}\|_2^2} \;.
\end{align*}}%
The inequality \eqref{dec13f} follows on choosing $k_1,\dots,k_n$ as
in \eqref{dec17a}.
\end{proof}

Lemma \ref{dec25b} does not apply directly to the unit ball, but it
has a counterpart which is applicable to the unit ball.

  \begin{lemma}\label{dec27b}
    Let $L\in\lhk$, where $\cH$ and $\cK$ are Hilbert spaces and $\cK
    = \cH \oplus \cdots \oplus \cH$ with~$d$ summands.  For each $j =
    1,\dots,d$, let $E_j \in\lkh$ be the operator $E_j (h_1
    \oplus\cdots\oplus h_d) = h_j$.  Suppose that $n = n_1 + \cdots +
    n_d$ for some nonnegative integers $n_1,\dots,n_d$.  Assume that
    $n \ge 2$, and set
    \begin{equation}\label{dec27c}
      K = \sum E_{j_1} L E_{j_2} L \cdots E_{j_{n-1}} L E_{j_n},
    \end{equation}
    where the summation is over all distinct arrangements
    $(j_1,j_2,\dots,j_n)$ of the tuple
    \begin{equation}\label{dec27d}
      (\,
    \overbrace{1,\dots,1}^{n_1}, \overbrace{2,\dots,2}^{n_2}, \dots,
    \overbrace{d,\dots,d}^{n_d}\, ) \;.
    \vphantom{\underbrace{x}{x}}
    \end{equation}
    Then
    \begin{equation}\label{dec27e}
          \| K \| \le d^{\,(n-1)/2}\, \| L \|^{n-1} \;.                     
    \end{equation}
  \end{lemma}

  \begin{proof}[Proof of Lemma~\ref{dec27b}]
    The argument is similar to the proof of Lemma~\ref{dec25b}.  We
    show how it works when $n_1=3$, $n_2=2$, and $d=2$.  Write the
    adjoint of $K$ as
    \begin{multline*}
      K^* = \tK[1,1,1,2,2] = 
       \sum E_{j_5}^* L^* E_{j_4}^* L^*E_{j_3}^* 
                          L^* E_{j_2}^*L^* E_{j_1}^* \\
          = \sum  E_{j_1}^* L^* E_{j_2}^* L^*E_{j_3}^* 
                          L^* E_{j_4}^*L^* E_{j_5}^* \;,
    \end{multline*}
    where the summations are over all distinct arrangements
    of~$(1,1,1,2,2)$.  Define the expressions \smash{$\tK[1,1,2,2],
      \tK[1,1,2], \dots$} in a similar way for the tuples
    $(1,1,2,2), (1,1,2), \dots$\;.  For any $f\in\cH$,
    \begin{equation*}
      K^*f = E_1^*L^*\tK[1,1,2,2]f + E_2^*L^*\tK[1,1,1,2] f .
    \end{equation*}
    Our assumptions imply that $E_1^*$ and $E_2^*$ have orthogonal
    ranges in $\cK$, and hence
    \begin{align*}
      \| K^*f \|^2 &= \| E_1^*L^*\tK[1,1,2,2]f \|^2
                           + \| E_2^*L^*\tK[1,1,1,2] f \|^2 \\[5pt]
&\le
      \| L^* \|^2 \Big\{ \| \tK[1,1,2,2]f \|^2
                           + \| \tK[1,1,1,2] f \|^2 \Big\} \\[5pt]
&=
    \| L^* \|^2 \Big\{ \| E_1^*L^*\tK[1,2,2]f \|^2
                           + \| E_2^*L^*\tK[1,1,2]f \|^2 \\[5pt]
&\hskip3cm
                     + \|E_1^*L^* \tK[1,1,2] f \|^2 
                           + \| E_2^*L^*\tK[1,1,1] f \|^2 \Big\} \;.
    \end{align*}
    By obvious estimates of the  separate terms and repetition of the
    process, we obtain
    \allowdisplaybreaks{
    \begin{align*}
      \| K^*f \|^2 
&\le 
    2\, \| L^* \|^4 \Big\{
      \| \tK[1,2,2]f \|^2
                           + \| \tK[1,1,2]f \|^2
            + \| \tK[1,1,1] f \|^2
    \Big\} \\
&=
    2\, \| L^* \|^4 \Big\{
      \|E_1^*L^*  \tK[2,2]f \|^2 +\|E_2^*L^* \tK[1,2]f \|^2 \\*
&\hskip2cm
    + \|E_1^*L^*  \tK[1,2]f \|^2 +\|E_2^*L^* \tK[1,1]f \|^2 
    + \|E_1^*L^*  \tK[1,1]f \|^2
    \Big\} \\
&\le
  2^2\, \| L^* \|^6 \Big\{
      \|\tK[2,2]f \|^2 +\|\tK[1,2]f \|^2 +\| \tK[1,1]f \|^2
    \Big\}  \\
&=
  2^2\, \| L^* \|^6 \Big\{
      \| E_2^*L^*E_2^*f \|^2 
      + \| E_1^*L^*E_2^*f \|^2  \\
&\hskip7cm + \| E_2^*L^*E_1^*f \|^2 
      +\| E_1^*L^*E_1^*f \|^2
    \Big\}  \\
&\le
   2^3\, \| L^* \|^8 \Big\{
    \| E_2^*f \|^2 + \| E_1^*f \|^2
    \Big\}  \\
&=
  2^4\, \| L^* \|^8 \| f \|^2.
    \end{align*}}%
Hence $\| K^* \| \le 2^2 \| L^* \|^4$, yielding \eqref{dec27e} when
$n_1=3$, $n_2=2$, and $d=2$.  

In the general case, we proceed in a similar way.  At the last stage
we shall have $\| K^*f \|^2 \le d^{n-1} \| L^* \|^{2(n-1)} \| f\|^2$,
which implies \eqref{dec27e}.
  \end{proof}

\begin{theorem}\label{ball}
  Let $\varphi(z)$ be a complex-valued function in the Schur-Agler
  class on the unit ball.  Then for each
  $z\in\bbB_d$ and any nonnegative integers $n_1,\dots,n_d$ which are
  not all zero,
  \begin{equation}\label{dec27f}
    \left|
      \frac{\partial^{\,n} \varphi}
    {\partial z_d^{n_d} \cdots \partial z_1^{n_1}}
    \right|
\le
    d^{\,(n-1)/2}\, n_1! n_2! \cdots n_d! \; \frac{1-|\varphi(z)|^2}
     {(1-\|z\|_2^2 )(1- \|z\|_2)^{n-1}} \,,
  \end{equation}
  where $n = n_1 + \cdots + n_d$\,.
\end{theorem}

This is an improvement on the inequality
\begin{equation*}
      \left|
      \frac{\partial^{\,n} \varphi}
    {\partial z_d^{n_d} \cdots \partial z_1^{n_1}}
    \right|
\le
    n!
    \; \frac{1-|\varphi(z)|^2}
     {(1-\|z\|_2^2 )(1- \|z\|_2)^{n-1}} \,,
\end{equation*}
which follows as a corollary of Theorem \ref{dec13d}.

  \begin{proof}[Proof of Theorem~\ref{ball}]
    When $n=1$, \eqref{dec27f} is immediate from \eqref{dec13f}.
    Assume $n \ge 2$.  Represent $\varphi(z)$ as in \eqref{aug23a}
    with the same choices of Hilbert spaces and operators as in the
    proof of Theorem~\ref{dec13d}.  Fix $z \in\bbB_d$, and set
    $L=L(z)$.  By \eqref{dec27a},
  \begin{equation*}
    \left|
      \frac{\partial^{\,n} \varphi}
    {\partial z_d^{n_d} \cdots \partial z_1^{n_1}}
    \right|
\le
  \| C (I_\cH - Z(z)A)^{-1}   \|  
  \;\| K \|  
  \;\|  (I_\cK -AZ(z))^{-1} B \|  ,
  \end{equation*}
  where $K$ is as in Lemma~\ref{dec27b}.  By \eqref{jun1f},
  \eqref{jun1g}, and Lemma~\ref{dec27b},
  \begin{equation*}
        \left|
      \frac{\partial^{\,n} \varphi}
    {\partial z_d^{n_d} \cdots \partial z_1^{n_1}}
    \right|
\le
  \frac{1 - |\varphi(z)|^2}{1 - \|Z(z)\|^2}\;
     d^{\,(n-1)/2}\, \| L \|^{n-1} .
  \end{equation*}
  Hence since $\| L \| \le 1/(1-\|Z(z)\|)$ by \eqref{dec26c} and
  $\|Z(z)\| = \|z\|_2$ by the proof of Theorem~\ref{dec13d}, we obtain
  the inequality \eqref{dec27f}.
  \end{proof}
  
  It is not hard to find cases where either one of the inequalities
  \eqref{dec13d} and \eqref{dec27f} is better than the other.  

\section{Concluding remarks}\label{S:remarks}

\subsection*{Operator-valued functions}
The results of \S\ref{S:polydisk} and \S\ref{S:ball} extend easily to
norm estimates of operator-valued functions which belong to the
Schur-Agler class on the polydisk or unit ball.  For the definitions
of these classes, see Agler \cite{agler1990} and Ball and Trent
\cite{balltrent} in the case of the polydisk, and Ball, Trent, and
Vinnikov \cite{btv2001} and Eschmeier and Putinar \cite{ep2002} for
the unit ball.  Operator-valued functions in the Schur-Agler class on
the polydisk or unit ball have unitary realizations of the same form
as described in \S\ref{polydiskcase} and \S\ref{ballcase} for the
scalar case.  Since such realizations are all that is needed for our
main estimates, the results of \S\ref{S:polydisk} and \S\ref{S:ball}
carry over without essential change to operator-valued functions.

\subsection*{Some open problems}
Many questions remain open, for example:

\noindent 1.
The estimates in \S\ref{S:polydisk} and \S\ref{S:ball} are best
possible in one variable \cite{AR,ruscheweyh}.  The corresponding
questions remain to be investigated in general.  We note that equality
holds at the origin in \eqref{dec20b} if $\varphi(z) = z_1^{n_1}
\cdots z_d^{n_d}$ for some nonnegative integers $n_1,\dots,n_d$ which
are not all zero.  This function belongs to the Schur-Agler class on
the polydisk by the definition of the class in \S\ref{polydiskcase}.

\noindent 2.
Are there functions in the Schur class (but not, of course, in the
Schur-Agler class) on the polydisk or unit ball that do not satisfy
the estimates in \S\ref{S:polydisk} and \S\ref{S:ball}?  A well-known
example of a function in the Schur class on $\bbD^3$ which is not in
the Schur-Agler class on the polydisk is the Kaijser-Varopoulos
polynomial \cite{lottosteger,varopoulos},
\begin{equation*}
  \varphi(z) = \frac15\; \big( z_1^2 + z_2^2 + z_3^2
                 - 2 z_1z_2 - 2 z_1 z_3 -2 z_2z_3   \big) \,.
\end{equation*}
We do not know if this function satisfies the inequalities in
\S\ref{S:polydisk}.  Examples of functions which belong to the Schur
class on the unit ball $\bbB_2$ but are not in the Schur-Agler class
on the unit ball are given by Alpay and Kaptano{\u{g}}lu
\cite[p.~6]{alpaykap}.  They have the form
\begin{equation*}
  \varphi(z) = z_1 + c_1 z_2^2 + c_2 z_2^4 + \cdots + c_m z_2^{2m} ,
\end{equation*}
where $m$ is an arbitrary positive integer and $c_1,c_2,\dots$ are the
coefficients in the expansion $1 - \sqrt{1-t} = c_1 t + c_2 t^2 +
\cdots\,$.  For each $j \ge 1$, $c_j > 0$.  We do not know if these
functions satisfy the inequalities in \S\ref{S:ball}.  Other simple
examples for the unit ball are given in \cite[pp.\ 
99--100]{aglermccarthy}.

\noindent 3.  The Schur-Agler class on the polydisk and unit ball are
two concrete classes of analytic functions which admit realizations
\eqref{may29a} such that the functions \eqref{jun2a} and \eqref{jun2b}
are linear.  Similar ingredients are available in Cartan domains, and
this environment may be amenable to analogous treatment.  See Upmeier
\cite{upmeier} for information on Cartan domains.  Realization
formulas were first proved for Cartan domains and their products by
Ambrozie and Timotin in \cite{AT}.  This was generalized by Ball and
Bolotnikov \cite{BB} to domains such that $1-Z(z)^*Z(z) > 0$, where
$Z(z)$ is a matrix polynomial.



\def\cprime{$'$} \def\cprime{$'$}
\providecommand{\bysame}{\leavevmode\hbox to3em{\hrulefill}\thinspace}
\providecommand{\MR}{\relax\ifhmode\unskip\space\fi MR }
\providecommand{\MRhref}[2]{%
  \href{http://www.ams.org/mathscinet-getitem?mr=#1}{#2}
}
\providecommand{\href}[2]{#2}

\bibliographystyle{amsplain}

\end{document}